\documentclass[12pt,reqno,a4wide]{article}
\usepackage{eurosym}
\usepackage{a4wide}
\usepackage{comment}
\usepackage{mathptmx}
\usepackage{amssymb}
\usepackage{amsfonts}
\usepackage{amsthm}
\usepackage{amsmath}
\usepackage{dsfont}
\usepackage{graphicx}
\usepackage{shadow}
\usepackage[all]{xy}
\usepackage{enumerate}
\usepackage{makeidx}
\usepackage{mathrsfs}
\usepackage{upgreek}
\usepackage{stmaryrd}
\usepackage[utf8]{inputenc}
\usepackage{array}
\usepackage{tikz-cd}
\usepackage{thumbpdf}
\usepackage[colorlinks=true,pagebackref]{hyperref}
\usepackage{textcomp}
\usepackage{chapterbib}
\usepackage{graphics}
\usepackage{longtable}
\usepackage{epsf}
\usepackage{bm}
\usepackage{adjustbox}
\usepackage{centernot}

\makeindex
\newtheorem{theorem}{Theorem}[section]
\newtheorem{lem}[theorem]{Lemma}
\newtheorem{thm}[theorem]{Theorem}
\newtheorem{prop}[theorem]{Proposition}
\newtheorem{cor}[theorem]{Corollary}
\theoremstyle{definition}
\newtheorem*{Beweis}{Proof}
\newtheorem{defn}[theorem]{Definition}
\newtheorem{definition}[theorem]{Definition}
\newtheorem{defns}[theorem]{Definitions}
\newtheorem{rem}[theorem]{Remark}
\newtheorem{rems}[theorem]{Remarks}
\newtheorem{punto}[theorem]{}

\theoremstyle{remark}

\newtheorem{ex}[theorem]{Example}
\newtheorem{exs}[theorem]{Examples}

\input{tcilatex}
\begin{document}

\title{On Semisimple Semirings\thanks{%
MSC2010: Primary 16D60; Secondary 16D40, 16D50 \newline
Key Words: Semirings; Semimodules; Semisimple Semirings; Injective
Semimodules; Projective Semimodules \newline
The authors would like to acknowledge the support provided by the Deanship
of Scientific Research (DSR) at King Fahd University of Petroleum $\&$
Minerals (KFUPM) for funding this work through projects No. RG1304-1 $\&$
RG1304-2}}
\author{$%
\begin{array}{ccc}
\text{Jawad Abuhlail}\thanks{\text{Corresponding Author}} &  & \text{Rangga
Ganzar Noegraha}\thanks{\text{The paper is extracted from his Ph.D.
dissertation under the supervision of Prof. Jawad Abuhlail.}} \\
\text{abuhlail@kfupm.edu.sa} &  & \text{rangga.gn@universitaspertamina.ac.id}
\\
\text{Department of Mathematics and Statistics} &  & \text{Universitas
Pertamina} \\
\text{King Fahd University of Petroleum $\&$ Minerals} &  & \text{Jl. Teuku
Nyak Arief} \\
\text{31261 Dhahran, KSA} &  & \text{Jakarta 12220, Indonesia}%
\end{array}%
$}
\date{\today }
\maketitle

\begin{abstract}
We investigate\emph{\ ideal-semisimple} and \emph{congruence-semisimple}
semirings. We give several new characterizations of such semirings using $e$%
-projective and $e$-injective semimodules. We extend several
characterizations of semisimple rings to (not necessarily subtractive)
commutative semirings.
\end{abstract}

%\addcontentsline{toc}{section}{\protect\numberline{}{Introduction}}

\section*{Introduction}

\emph{Semirings} (defined, roughly, as rings not necessarily with
subtraction) can be considered as a generalization of both rings and
distributive bounded lattices. Semirings{, and their \emph{semimodules}
(defined, roughly, as modules not necessarily with subtraction), have wide
applications in many aspects of Computer Science and Mathematics, e.g.,
Automata Theory \cite{HW1998}, Tropical Geometry \cite{Gla2002} and
Idempotent Analysis \cite{LM2005}. Many of these applications can be found
in Golan's book \cite{Gol1999}, which is considered a main reference in this
topic. }

\bigskip

Several papers by Abuhlail, I'llin, Katsov and Nam (among others) prepared
the stage for a homological characterization of special classes of semirings
using special classes of projective, injective and flat semimodules (cf.,
\cite{KNT2009}, \cite{Ili2010}, \cite{KN2011}, \cite{Abu2014}, \cite{KNZ2014}%
, \cite{AIKN2015}, \cite{IKN2017}, \cite{AIKN2018}).

\bigskip

The notions of projective and injective objects can be defined in any
category relative to a suitable \emph{factorization system} of its arrows.
Projective, injective and flat semimodules have been studied intensively
(see \cite{Gla2002} for details). Recently, left (right) $V$-semirings, all
of whose \emph{congruence-simple} left (right) semimodules are injective
have been completely characterized in \cite{AIKN2015}, and \emph{%
ideal-semisimple} semirings all of whose left cyclic semimodules are
projective have been investigated in \cite{IKN2017}.

\bigskip

In addition to the \emph{categorical notions} of \emph{projective} and \emph{%
injective} \emph{semimodules} over a semiring, new notions of projectivity
and injectivity of semimodules over semirings were considered by the first
author who introduced the so called $e$\emph{-projective\ }and $e$\emph{%
-injective semimodules} \cite{Abu2014-CA}. One reason for the interest in
such notions is the phenomenon that assuming that \emph{all} semimodules of
a given semiring $S$ are projective (injective) forces the semiring to be a
\emph{ring} (cf., \cite[Theorem 3.4]{Ili2010}).

\bigskip

The paper is divided into two sections.

\bigskip

In Section 1, we collect the basic definitions, examples and preliminaries
used in this paper. Among others, we include the definitions and basic
properties of \emph{exact sequences }as defined by Abuhlail \cite{Abu2014}.

\bigskip

In Section 2, we investigate \emph{ideal-semisimple} and \emph{%
congruence-semisimple} semirings. A semiring $S$ is left (right)
ideal-semisimple, iff $S$ is a direct sum of ideal-simple left (right)
ideals. By \cite[Theorem 7.8]{HW1996}, $S$ is left (right) ideal-semisimple
if and only if $S\simeq M_{n_{1}}(D_{1})\times \cdots \times M_{n_{k}}(D_{k})
$, where $D_{i}$ is a division semiring and $M_{n_{i}}(D_{i})$ is the
semiring of $n_{i}\times n_{i}$-matrices over $D_{i}$ for each $i=1,\cdots
,k.$ A \emph{left subtractive} semiring $S$ was shown to be left
ideal-semisimple if and only if every left $S$-semimodule is $S$-$k$%
-projective ($S$-$k$-injective) \cite[Theorem 4.4]{KNT2009}. In Proposition %
\ref{proj-impl} (Proposition \ref{sum-einj}), we show that a semiring $S$
over which every left $S$-semimodule is $S$-$k$-projective ($S$-$k$%
-injective) is a finite direct sum of irreducible summands.

\bigskip

In Section 3, we restrict our attention to commutative semirings. In Theorem %
\ref{ISS-COMM}, we extend several classical characterizations of semisimple
semirings to commutative, \emph{not necessarily subtractive}, semirings. In
Theorem \ref{CSS-COMM}, we show that a \emph{commutative} semiring $S$ is
ideal-semisimple if and only if every $S$-semimodule is $S$-$e$-injective ($%
S $-$k$-injective) and $S$ satisfies some technical condition. The two
results are combined in Theorem \ref{com-idss} to provide a complete
characterization of commutative ideal-semisimple semirings. The \emph{%
congruence-semisimple} version of this main result is given in Theorem \ref%
{com-css}. Examples \ref{B(3,1)} and \ref{nonseproj} demonstrate that the
conditions assumed in our main results in this section, in particular the
commutativity of the base semiring, cannot be dropped.

\section{Preliminaries}

\label{prelim}

\qquad In this section, we provide the basic definitions and preliminaries
used in this work. Any notions on semirings and semimodules that are not
defined can be found in {our main reference \cite{Gol1999}. We refer to \cite%
{Wis1991} for the foundations of Module and Ring Theory.}

\begin{definition}
(\cite{Gol1999}) A \textbf{semiring}%
\index{Semiring} is a datum $(S,+,0,\cdot ,1)$ consisting of a commutative
monoid $(S,+,0)$ and a monoid $(S,\cdot ,1)$ such that $0\neq 1$ and%
\begin{eqnarray*}
a\cdot 0 &=&0=0\cdot a%
\text{ for all }a\in S; \\
a(b+c) &=&ab+ac\text{ and }(a+b)c=ac+bc\text{ for all }a,b,c\in S.
\end{eqnarray*}
\end{definition}

\begin{punto}
\cite{Gol1999} Let $S$ and $T$ be semirings. The categories $_{S}\mathbf{SM}$
of \textbf{left} $S$-\textbf{semimodules} with arrows the $S$-linear maps, $%
\mathbf{SM}_{T}$ of right $S$-semimodules with arrows the $T$-linear maps,
and $_{S}\mathbf{SM}_{T}$ of $(S,T)$-bisemimodules are defined in the usual
way (as for modules and bimodules over rings). We write $L\leq _{S}M$ to
indicate that $L$ is a subsemimodule of the left (right) $S$-semimodule $M.$
\end{punto}

\begin{defns}
\label{def-semiring}(\cite{Gol1999}) Let $(S,+,0,\cdot ,1)$ be a semiring.

\begin{itemize}
\item If the monoid $(S,\cdot ,1)$ is commutative, we say that $S$ is a
\emph{commutative semiring}. If moreover, $(S\backslash \{0\},\cdot ,1)$ is
a group, we say that $S$ is a \textbf{semifield}.

\item Let%
\begin{equation}
V(S):=\{s\in S\mid s+t=0\text{ for some }t\in S\}.  \label{V(S)}
\end{equation}%
If $V(S)=\{0\}$, we say that $S$ is \textbf{zerosumfree}. Notice that $%
V(S)=S $ if and only if $S$ is a ring.

\item The set of \emph{cancellative elements of }$S$ is defined as%
\begin{equation*}
K^{+}(S)=\{x\in S\mid x+y=x+z\Longrightarrow y=z\text{ for any }y,z\in S\}.
\end{equation*}%
We say that $S$ is a \textbf{cancellative semiring,} iff $K^{+}(S)=S.$
\end{itemize}
\end{defns}

\begin{exs}
({\cite{Gol1999})}

\begin{itemize}
\item Every ring is a cancellative semiring.

\item Any \emph{distributive bounded lattice} $\mathcal{L}=(L,\vee ,\wedge
,0,1)$ is a commutative semiring.

\item $(\mathbb{Z}^{+},+,0,\cdot ,1)$ (resp. $(\mathbb{Q}^{+},+,0,\cdot ,1),$
$(\mathbb{R}^{+},+,0,\cdot ,1)$), the set of non-negative integers (resp.
non-negative rational numbers, non-negative real numbers) is a commutative
cancellative semiring which is not a ring.

\item $M_{n}(S),$ the set of all $n\times n$ matrices over a semiring $S,$
is a semiring.

\item $\mathbb{B}:=\{0,1\}$,with $1+1=1,$ is a semiring called the \textbf{%
Boolean semiring}. $\mathbb{B}$ is an semifield which is not a field.

\item The \emph{max-plus algebra} $\mathbb{R}_{\max ,+}:=(\mathbb{R}\cup
\{-\infty \},\max ,-\infty ,+,0)$ is an additively idempotent semiring.

\item The \emph{log algebra} $(\mathbb{R}\cup \{-\infty ,\infty \},\oplus
,\infty ,+,0)$ is a semiring, where%
\begin{equation*}
x\oplus y=-ln(e^{-x}+e^{-y})
\end{equation*}
\end{itemize}
\end{exs}

\begin{ex}
\label{B(n,i)}(\cite[Example 1.8]{Gol1999}, \cite{AA1994}) Consider
\begin{equation*}
B(n,i):=(B(n,i),\oplus ,0,\odot ,1),
\end{equation*}%
where $B(n,i)=\{0,1,2,...,n-1\}$ and%
%($n\geq 2$) is considered
%with the binary operations $\oplus $ and $\odot $ are defined as follows:

$a\oplus b=a+b$ if $a+b<n$; otherwise, $a\oplus b=c$ with $i\leq c<n$ is the
unique natural number satisfying $c\equiv (a+b)$ mod $(n-i);$

$a\odot b=ab$ if $ab<n$; otherwise, $a\odot b=c$ with $i\leq c<n$ is the
unique natural number with $c\equiv ab$ mod $(n-i)$.

Then $B(n,i)$ is a semiring. Notice that $B(n,0)=\mathbb{Z}_{n}$ (a group)
and that $B(2,1)=\mathbb{B}$ (the Boolean Algebra).
\end{ex}

\begin{ex}
(\cite[page 150, 154]{Gol1999}) Let $S$ be a semiring, $M$ be a left $S$%
-semimodule and $L\subseteq M.$ The \textbf{subtractive closure }of $L$ is
defined as%
\begin{equation}
\overline{L}:=\{m\in M\mid \text{ }m+l=l^{\prime }\text{ for some }%
l,l^{\prime }\in L\}.  \label{L-s-closure}
\end{equation}%
One can easily check that $\overline{L}=Ker(M\overset{\pi }{\longrightarrow }%
M/L),$ where $\pi $ is the canonical projection. We say that $L$ is \textbf{%
subtractive}, iff $L=\overline{L}.$ The left $S$-semimodule $M$ is a \textbf{%
subtractive semimodule}, iff every $S$-subsemimodule $L\leq _{S}M$ is
subtractive.
\end{ex}

\begin{definition}
\cite[page 71]{Gol1999} Let $S$ be a semiring. We say that $S$ is a \textbf{%
left subtractive semiring }(\textbf{right subtractive semiring}), iff every
left (right) ideal of $S$ is subtractive. We say that $S$ is a \textbf{%
subtractive semiring}, iff $S$ is both left and right subtractive.
%On the other hand, we
%say that $S$ is a \textbf{left }(a \textbf{right austere semiring}) if $%
%\{0\}\ $is the only subtractive left (right) ideal of $S.$ We say that $S$
%is an \textbf{austere semiring} if $S$ is left austere and right austere.
\end{definition}

\begin{rem}
Whether a left subtractive semiring is necessarily right subtractive was an
open problem till a counterexample was given in \cite[Fact 2.1]{KNT2011}.
\end{rem}

Following \cite{BHJK2001}, we use the following definitions.

\begin{punto}
\label{variety}(cf., \cite{AHS2004})\ The category $_{S}\mathbf{SM}$ of left
semimodules over a semiring $S$ is a \emph{variety} in the sense of
Universal Algebra (closed under homomorphic images, subobjects and arbitrary
products). Whence $_{S}\mathbf{SM}$ is complete, i.e. has all limits (e.g.,
direct products, equalizers, kernels, pullbacks, inverse limits) and
cocomplete, i.e. has all colimits (e.g., direct coproducts, coequalizers,
cokernels, pushouts, direct colimits).
\end{punto}

\subsection*{Semisimple Semimodules}

\begin{punto}
\cite[page 162]{Gol1999} Let $S$ be a semiring.

An equivalence relation $\rho $ on a left $S$-semimodule $M$ is a \textbf{%
congruence relation}, iff it preserves the addition and the scalar
multiplication on $M,$ \emph{i.e. }for all $s\in S$ and $m,m^{\prime
},n,n^{\prime }\in M:$%
\begin{equation*}
m\rho m^{\prime }\text{ and }n\rho n^{\prime }\Longrightarrow (m+m^{\prime
})\rho (n+n^{\prime }),
\end{equation*}%
\begin{equation*}
m\rho m^{\prime }\Longrightarrow (sm)\rho (sm^{\prime }).
\end{equation*}%
A \textbf{congruence relation on the semiring }$S$ is an equivalence
relation $\rho $ on $S$ such that or all $s,s^{\prime },t,t^{\prime }\in S:$%
\begin{equation*}
s\rho t\text{ and }s^{\prime }\rho t^{\prime }\Longrightarrow (s+s^{\prime
})\rho (t+t^{\prime })\text{ and }(ss^{\prime })\rho (tt^{\prime }).
\end{equation*}
\end{punto}

\begin{ex}
Let $S$ be a semiring, $M$ a left $S$-semimodule and $N\leq _{S}M.$ The
\textbf{Bourne relation} $\equiv _{N}$ on $M$ is defined as:
\begin{equation*}
m\equiv _{N}m^{\prime }\Leftrightarrow m+n=m^{\prime }+n^{\prime }\text{ for
some }n,n^{\prime }\in N.
\end{equation*}%
It is clear that $\equiv _{N}$ is a congruence relation. Moreover, $%
M/N=M/\equiv _{N}=\{[m]_{N}$ $|$ $m\in M\}$ $(=M/\overline{N})$ is a left $S$%
-semimodule, the canonical surjective map $\pi _{N}:M\longrightarrow M/N$ is
$S$-linear, and $Ker(\pi _{N})=\overline{N}.$ In particular, $Ker(\pi
_{N})=0 $ if and only if $N\leq _{S}M$ is subtractive (this explains why
subtractive ideals are called $k$-ideals in many references).
\end{ex}

Following \cite{BHJK2001}, we use the following definitions.

\begin{definition}
Let $S$ be a semiring. A left $S$-semimodule $M$ is

\textbf{ideal}-\textbf{simple}, iff 0 and $M$ are the only $S$%
-subsemimodules of $M$;

\textbf{congruence-simple}, iff
\begin{equation*}
\Delta _{M}:=\{(m,m)\mid m\in M\}
\end{equation*}
and $M\times M$ are the only congruence relations on $M$.
\end{definition}

\begin{defn}
We say that the semiring $S$ is

\textbf{left ideal-simple }(\textbf{right ideal-simple}), iff $0$ and $S$
are the only left (right) ideals of $S;$ equivalently, $S$ is ideal-simple
as a left (right) $S$-semimodule.;

\textbf{left congruence-simple} (\textbf{right congruence-simple}), iff $S$
is congruence-simple as a left (right) $S$-semimodule;

\textbf{ideal-simple}, iff $0$ and $S$ are the only (two-sided) ideals of $%
S; $

\textbf{congruence-simple}, iff $\Delta _{S}$ and $S\times S$ are the only
congruence relations on the semiring $S.$
\end{defn}

\begin{rem}
If $M$ is a congruence-simple left $S$-semimodule, then the only subtractive
$S$-subsemimodules of $M$ are $0$ and $M.$ To show this, suppose that $N\neq
0$ is a subtractive $S$-subsemimodule of $M.$ Then $\equiv _{N}$ is a
congruence relation on $M$ with $n\equiv _{N}0$ for some $n\in N\backslash 0$%
. Thus $\equiv _{N}\neq \Delta _{M}$, which implies $\equiv _{N}=M^{2}$ as $%
M $ is congruence-simple. If $m\in M$, then $mM^{2}0$, that is $m\equiv
_{N}0 $. Therefore, there exist $n,n^{\prime }\in N$ such that $%
m+n=n^{\prime }$. Since $N$ subtractive, $m\in N$. Hence $M=N$.
\end{rem}

\begin{ex}
\label{cong-not-id}\cite[3.7 (b)]{KNZ2014} Let $(M,+,0)$ be a finite lattice
that is not distributive. The endomorphism semiring $\mathbf{E}_{M}$ of $M$
is a congruence-simple semiring which is \emph{not} ideal-simple.
\end{ex}

\begin{ex}
\label{id-not-cong}(cf., \cite[3.7 (c)]{KNZ2014})\ Every zerosumfree
division semiring $D$ that is not isomorphic to $\mathbb{B}$ (e.g., $\mathbb{%
R}^{+}$) is an ideal-simple semiring which is \emph{not} congruence-simple
since%
\begin{equation*}
\rho =\{(a,b)|\text{ }a,b\in D\backslash \{0\}\}\cup \{(0,0)\}
\end{equation*}%
is a non-trivial non-universal congruence relation on $D.$
\end{ex}

\begin{lem}
\label{cong-s-char}A left $S$-semimodule $M$ is congruence-simple if and
only if every non-zero $S$-linear map from $M$ is injective.
\end{lem}

\begin{Beweis}
($\Rightarrow $) Let $f:M\rightarrow N$ be a non-zero $S$-linear map and
consider the congruence relation on $M$ defined by%
\begin{equation*}
m\equiv _{f}m^{\prime }\Longleftrightarrow f(m)=f(m^{\prime }).
\end{equation*}%
Pick some $m\in M\backslash \{0\}$ such that $f(m)\neq 0$. Since $\equiv
_{f} $ is a congruence relation on $M$ with $m\not\equiv _{f}0$, we know $%
\equiv _{f}\neq M^{2}.$ It follows that $\equiv _{f}=\Delta _{M}$ as $M$ is
congruence-simple. Hence $f$ is injective.

($\Leftarrow $) Assume that $M$ is congruence-simple. Let $\rho $ be a
congruence relation on $M$. The canonical map $f:M\rightarrow M/\rho $ is $S$%
-linear. If $f=0$, then $[m]_{\rho }=[0]_{\rho }$ for every $m\in M$, that
is $m\rho 0$ for every $m\in M$ and $m\rho m^{\prime }$ for every $%
m,m^{\prime }\in M$. If $f\neq 0$, then $f$ is injective, that is $[m]_{\rho
}\neq \lbrack m^{\prime }]_{\rho }$ whenever $m\neq m^{\prime }$. Thus $%
m\rho m^{\prime }$ whenever $m\neq m^{\prime }$ and $\rho =\Delta _{M}$.$%
\blacksquare $
\end{Beweis}

\begin{lem}
\label{id-ss-char}A left $S$-semimodule $M$ is ideal-simple if and only if
every non-zero $S$-linear map to $M$ is surjective.
\end{lem}

\begin{Beweis}
($\Rightarrow $) Let $f:L\rightarrow M$ be a non-zero $S$-linear map. Then
there exists $l\in L\backslash 0$ such that $f(l)\neq 0$. Thus, $f(L)$ is a
non-zero subsemimodule of $M$ and so $f(L)=M$ as $M$ ideal-simple.

($\Leftarrow $) Let $K$ be a subsemimodule of $M$. Then the embedding $%
f:K\rightarrow M$ is an $S$-linear map. If $f=0$, then $K=f(K)=0$. If $f\neq
0$, then $f$ is surjective, that is $K=f(K)=M$.$\blacksquare $
\end{Beweis}

\begin{punto}
(\cite[page 184]{Gol1999}) Let $S$ be a semiring. A left $S$-semimodule $M$
is the \textbf{direct sum} of a family $\{L_{\lambda }\}_{\lambda \in
\Lambda }$ of $S$-subsemimodules $L_{\lambda }\leq _{S}M,$ and we write $%
M=\bigoplus\limits_{\lambda \in \Lambda }L_{\lambda }$, iff every $m\in M$
can be written in a \emph{unique way} as a finite sum $m=l_{\lambda
_{1}}+\cdots +l_{\lambda _{k}}$ where $l_{\lambda _{i}}\in L_{\lambda _{i}}$
for each $i=1,\cdots ,k.$ Equivalently, $M=\bigoplus\limits_{\lambda \in
\Lambda }L_{\lambda }$ if $M=\sum\limits_{\lambda \in \Lambda }L_{\lambda }$
and for each finite subset $A\subseteq \Lambda $ with $l_{a},l_{a}^{\prime
}\in L_{a},$ we have:%
\begin{equation*}
\sum\limits_{a\in A}l_{a}=\sum\limits_{a\in A}l_{a}^{\prime }\Longrightarrow
l_{a}=l_{a}^{\prime }\text{ for all }a\in A.
\end{equation*}
\end{punto}

\begin{punto}
An $S$-semimodule $N$ is a \textbf{retract} of an $S$-semimodule $M,$ iff
there exists a (surjective) $S$-linear map $\theta :M\longrightarrow N$ and
an (injective) $S$-linear map $\psi :N\longrightarrow M$ such that $\theta
\circ \psi =\mathrm{id}_{N}$ (equivalently, $N\simeq \alpha (M)$ for some
\emph{idempotent} endomorphism $\alpha \in \mathrm{End}(M_{S})$).
\end{punto}

\begin{punto}
An $S$-semimodule $N$ is a \textbf{direct summand} of an $S$-semimodule $M$ (%
\emph{i.e.} $M=N\oplus N^{\prime }$ for some $S$-subsemimodule $N^{\prime }$
of $M$) if and only if there exists $\alpha \in \mathrm{Comp}(\mathrm{End}%
(M_{S}))$ s.t. $\alpha (M)=N$ where for any semiring $T$ we set%
\begin{equation*}
\mathrm{Comp}(T)=\{t\in T\mid \text{ }\exists \text{ }\widetilde{t}\in T%
\text{ with }t+\widetilde{t}=1_{T}\text{ and }t\widetilde{t}=0_{T}=%
\widetilde{t}t\}.
\end{equation*}%
Indeed, every direct summand of $M$ is a retract of $M;$ the converse is not
true in general; for example $N_{1}=\left\{ \left[ {%
\begin{array}{cc}
a & a \\
b & b%
\end{array}%
}\right] |\text{ }a,b\in \mathbb{R}^{+}\right\} $ is a retract of $M_{2}(%
\mathbb{R^{+}})$ that is not a direct summand. Golan \cite[Proposition 16.6]%
{Gol1999} provided characterizations of direct summands.
\end{punto}

\begin{rems}
\label{d-iso}Let $M$ be a left $S$-semimodule and $K,L\leq _{S}M$ be $S$%
-semimodules of $M.$

\begin{enumerate}
\item If $K+L$ is direct, then $K\cap L=0.$ The converse is not true in
general.

\item If $M=K\oplus L$, then $M/K\simeq L.$
\end{enumerate}
\end{rems}

\begin{ex}
Let $S=M_{2}(\mathbb{R}^{+})$. Notice that%
\begin{equation*}
E_{1}=\left\{ \left[ {%
\begin{array}{cc}
a & 0 \\
b & 0%
\end{array}%
}\right] |\text{ }a,b\in \mathbb{R}^{+}\right\} \text{ and }N_{\geq
1}\left\{ \left[ {%
\begin{array}{cc}
a & c \\
b & b%
\end{array}%
}\right] |\text{ }a\leq c,b\leq d,a,b,c,d\in \mathbb{R}^{+}\right\}
\end{equation*}%
are left ideals of $S$ with $E_{1}\cap N_{\geq 1}=\{0\}$. However, the sum $%
E_{1}+N_{\geq 1}$ is not direct since
\begin{equation*}
\left[ {%
\begin{array}{cc}
1 & 0 \\
0 & 0%
\end{array}%
}\right] +\left[ {%
\begin{array}{cc}
0 & 1 \\
0 & 0%
\end{array}%
}\right] =\left[ {%
\begin{array}{cc}
0 & 0 \\
0 & 0%
\end{array}%
}\right] +\left[ {%
\begin{array}{cc}
1 & 1 \\
0 & 0%
\end{array}%
}\right] .
\end{equation*}
\end{ex}

\begin{lem}
(\cite[Proposition 16.6]{Gol1999}) Let $S$ be a semiring. The following are
equivalent for a left $S$-semimodule $M$ and $N\leq _{S}M:$

\begin{enumerate}
\item $N$ is a direct summand of $M$ (i.e. $M=N\oplus N^{\prime }$ for some $%
S$-subsemimodule $N^{\prime }$ of $M$);

\item $N=\alpha (M)$ for some $\alpha \in \mathrm{Comp}(\mathrm{End}%
(M_{S})); $

\item $\exists $ $N^{\prime }\leq _{S}M$ such that $M=N+N^{\prime }$ and the
restrictions to $\equiv _{N}$ to $N^{\prime }$ and the restriction $\equiv
_{N^{\prime }}$ to $N$ are trivial.
\end{enumerate}
\end{lem}

\subsection*{Exact Sequences}

\bigskip

Throughout, $(S,+,0,\cdot ,1)$ is a semiring and, unless otherwise
explicitly mentioned, an $S$-module is a \emph{left }$S$-semimodule.

\bigskip

\begin{definition}
A morphism of left $S$-semimodules $f:L\rightarrow M$ is

$k$-\textbf{normal}, iff whenever $f(m)=f(m^{\prime })$ for some $%
m,m^{\prime }\in M,$ we have $m+k=m^{\prime }+k^{\prime }$ for some $%
k,k^{\prime }\in Ker(f);$

$i$-\textbf{normal}, iff $\func{Im}(f)=\overline{f(L)}$ ($:=\{m\in M|\text{ }%
m+l\in L\text{ for some }l\in L\}$).

\textbf{normal}, iff $f$ is both $k$-normal and $i$-normal.
\end{definition}

There are several notions of exactness for sequences of semimodules. In this
paper, we use the relatively new notion introduced by Abuhlail:

\begin{definition}
\label{Abu-exs}(\cite[2.4]{Abu2014}) A sequence
\begin{equation}
L\overset{f}{\longrightarrow }M\overset{g}{\longrightarrow }N  \label{LMN}
\end{equation}%
of left $S$-semimodules is

\textbf{exact}, iff $f(L)=Ker(g)$ and $g$ is $k$-normal;

\textbf{semi-exact}, iff $\overline{f(L)}=Ker(g);$

\textbf{proper-exact}, iff $f(L)=Ker(g).$
\end{definition}

\begin{punto}
We call a (possibly infinite) sequence of $S$-semimodules
\begin{equation}
\cdots \rightarrow M_{i-1}\overset{f_{i-1}}{\rightarrow }M_{i}\overset{f_{i}}%
{\rightarrow }M_{i+1}\overset{f_{i+1}}{\rightarrow }M_{i+2}\rightarrow \cdots
\label{chain}
\end{equation}

\emph{chain complex} if $f_{j+1}\circ f_{j}=0$ for every $j;$

\emph{exact} (resp., \emph{proper-exact}, \emph{semi-exact, quasi-exact}) if
each partial sequence with three terms $M_{j}\overset{f_{j}}{\rightarrow }%
M_{j+1}\overset{f_{j+1}}{\rightarrow }M_{j+2}$ is exact (resp.,
proper-exact, semi-exact, quasi-exact).

A \textbf{short exact sequence}%
\index{short exact sequence} (or a \textbf{Takahashi extension}%
\index{Takahashi extension} \cite{Tak1982b}) of $S$-semimodules is an exact
sequence of the form%
\begin{equation*}
0\longrightarrow L\overset{f}{\longrightarrow }M\overset{g}{\longrightarrow }%
N\longrightarrow 0
\end{equation*}%
The following result shows some of the advantages of the Abuhlail's
definition of exact sequences over the previous ones:
\end{punto}

\begin{lem}
\label{exact}Let $L,M$ and $N$ be $S$-semimodules.

\begin{enumerate}
\item $0\longrightarrow L\overset{f}{\longrightarrow }M$ is exact if and
only if $f$ is injective.

\item $M\overset{g}{\longrightarrow }N\longrightarrow 0$ is exact if and
only if $g$ is surjective.

\item $0\longrightarrow L\overset{f}{\longrightarrow }M\overset{g}{%
\longrightarrow }N$ is semi-exact and $f$ is normal (proper-exact and $f$ is
normal) if and only if $L\simeq \mathrm{Ker}(g).$

\item $0\longrightarrow L\overset{f}{\longrightarrow }M\overset{g}{%
\longrightarrow }N$ is exact if and only if $L\simeq \mathrm{Ker}(g)$ and $g$
is $k$-normal.

\item $L\overset{f}{\longrightarrow }M\overset{g}{\longrightarrow }%
N\longrightarrow 0$ is semi-exact and $g$ is normal if and only if $N\simeq
M/f(L).$

\item $L\overset{f}{\longrightarrow }M\overset{g}{\longrightarrow }%
N\longrightarrow 0$ is exact if and only if $N\simeq M/f(L)$ and $f$ is $i$%
-normal.

\item $0\longrightarrow L\overset{f}{\longrightarrow }M\overset{g}{%
\longrightarrow }N\longrightarrow 0$ is exact if and only if $L\simeq
\mathrm{Ker}(g)$ and $N\simeq M/L.$
\end{enumerate}
\end{lem}

\begin{cor}
\label{M/L}The following assertions are equivalent:

\begin{enumerate}
\item $0\rightarrow L\overset{f}{\rightarrow }M\overset{g}{\rightarrow }%
N\rightarrow 0$ is an exact sequence of $S$-semimodules;

\item $L\simeq \mathrm{Ker}(g)$ and $N\simeq M/f(L)$;

\item $f$ is injective, $f(L)=\mathrm{Ker}(g),$ $g$ is surjective and ($k$%
-)normal.

In this case, $f$ and $g$ are normal morphisms.
\end{enumerate}
\end{cor}

\section{Semisimple Semirings}

\label{sec-ss}

\markright{\scriptsize\tt Chapter \ref{sec-ss}: Semisimple
Semimodules and Semirings} %--------------------------------

Throughout, $(S,+,0,\cdot ,1)$ is a semiring and, unless otherwise
explicitly mentioned, an $S$-module is a \textbf{left} $S$-semimodule.

%The equivalence is a corollary of remark \ref{rem1}.

\begin{definition}
Let $S$ be a semiring. A left $S$-semimodule $M$ is called

\textbf{ideal}-\textbf{semisimple}, iff $M=\dbigoplus\limits_{\lambda \in
\Lambda }M_{\lambda },$ a direct sum of ideal-simple $S$-subsemimodules;

\textbf{congruence-semisimple}, iff $M=\dbigoplus\limits_{\lambda \in
\Lambda }M_{\lambda },$ a direct sum of congruence-simple $S$-subsemimodules.
\end{definition}

\begin{defn}
A semiring $S$ is

\textbf{left ideal-semisimple} (resp., \textbf{right ideal-semisimple}), iff
$S$ is ideal-semisimple as a left (right) $S$-semimodule, equivalently $S$
is a \emph{finite} direct sum of ideal-simple left (right) ideals.

\textbf{left congruence-semisimple} (resp., \textbf{right
congruence-semisimple}), iff $S$ is congruence-semisimple as a left (right) $%
S$-semimodule; equivalently $S$ is a \emph{finite} direct sum of
congruence-simple left (right) ideals.
\end{defn}

\begin{definition}
(\cite{AIKN2018})\ A left $S$-semimodules $P$ is

$M$-$e$-\textbf{projective} (where $M$ is a left $S$-semimodule), iff the
covariant functor
\begin{equation*}
Hom_{S}(P,-):%
\text{ }_{S}\mathbf{SM}\longrightarrow \text{ }_{\mathbb{Z}^{+}}\mathbf{SM}
\end{equation*}%
transfers every short exact sequence of left $S$-semimodules%
\begin{equation}
0\longrightarrow L\overset{f}{\longrightarrow }M\overset{g}{\longrightarrow }%
N\longrightarrow 0  \label{lmn-k}
\end{equation}%
into a short exact sequence of commutative monoids%
\begin{equation}
0\longrightarrow Hom_{S}(P,L)\overset{(P,f)}{\longrightarrow }Hom_{S}(P,M)%
\overset{(P,g)}{\longrightarrow }Hom_{S}(P,N)\longrightarrow 0.
\label{hom(p,lmn)}
\end{equation}%
We say that $P$ is $e$\textbf{-projective}, iff $P$ is $M$-$e$-projective
for every left $S$-semimodule $M$.
\end{definition}

\begin{punto}
Let $P$ be a left $S$-semimodule. For a left $S$-semimodule $M,$ we say that
$P$ is $M$\textbf{-projective} (resp. $M$-$k$\textbf{-projective}), iff for
every (\emph{normal}) surjective $S$-linear map $f:M\rightarrow N$ and an $S$%
-linear map $g:P\rightarrow N,$ there exists an $S$-linear map $%
h:P\rightarrow M$ such that $f\circ h=g.$

We say that $P$ is \textbf{projective} (resp., $k$\textbf{-projective}), iff
$P$ is $M$-projective (resp., $M$-$k$-projective) for every left $S$%
-semimodule $M.$
\end{punto}

\begin{definition}
(\cite[1.24]{Abu2014-CA}) Let $M$ be a left $S$-semimodule. A left $S$%
-semimodules $J$ is $M$-$e$-\textbf{injective}, iff the contravariant functor%
\begin{equation*}
Hom_{S}(-,J):\text{ }_{S}\mathbf{SM}\longrightarrow \text{ }_{\mathbb{Z}^{+}}%
\mathbf{SM}
\end{equation*}%
transfers every short exact sequence of left $S$-semimodules
\begin{equation*}
0\longrightarrow L\overset{f}{\longrightarrow }M\overset{g}{\longrightarrow }%
N\longrightarrow 0
\end{equation*}%
into a short exact sequence of commutative monoids%
\begin{equation*}
0\longrightarrow Hom_{S}(N,J)\longrightarrow Hom_{S}(M,J)\longrightarrow
Hom_{S}(L,J)\longrightarrow 0.
\end{equation*}%
We say that $J$ is $e$\textbf{-injective,} iff $J$ is $M$-$e$-injective for
every left $S$-semimodule $M.$
\end{definition}

\begin{punto}
Let $I$ be a left $S$-semimodule. For a left $S$-semimodule $M,$ we say that
$I$ is $M$\textbf{-injective} ($M$-$i$\textbf{-injective}) \cite[page 197]%
{Gol1999}, iff for every (\emph{normal}) injective $S$-linear map $%
f:L\rightarrow M$ and any $S$-linear map $g:L\rightarrow I,$ there exists an
$S$-linear map $h:M\rightarrow I$ such that $h\circ f=g.$

We say that $I$ is \textbf{injective} (resp., $i$\textbf{-injective}) if $I$
is $M$-injective (resp., $M$-$i$-injective) for every left $S$-semimodule $%
M. $
\end{punto}

{The following characterizations of \emph{semisimple rings} are well known
(cf., \cite{Gri2007}):}

\begin{thm}
\label{semi-simple-rings}\emph{(\cite[page 362, 402, 404]{Gri2007})} Let $R$
be a ring. Then the following assertions are equivalent:

\begin{enumerate}
\item Every left (right) $R$-module is $R$-injective;

\item Every left (right) $R$-module is injective;

\item Every left (right) $R$-module is projective;

\item Every short exact sequence of left (right) $R$-modules $0\rightarrow
L\rightarrow M\rightarrow N\rightarrow 0$ splits;

\item Every left (right) ideal of $R$ is a direct summand;

\item $R$ is left (right) semisimple.
\end{enumerate}
\end{thm}

\begin{lem}
\label{sumproj}\emph{(\cite[Lemma 3.15.]{AN-1})\ }If $M$ is a left $S$%
-semimodule such that every subtractive subsemimodule is a direct summand,
then every left $S$-semimodule is $M$-$e$-projective.
\end{lem}

\begin{punto}
We say that a sequence of $S$-semimodules%
\begin{equation}
0\rightarrow A\overset{f}{\longrightarrow }B\overset{g}{\longrightarrow }%
C\rightarrow 0  \label{se}
\end{equation}%
is

\emph{left splitting,} iff there exists $f^{\prime }\in \mathrm{Hom}%
_{S}(B,A) $ such that $f^{\prime }\circ f=id_{A};$

\emph{right splitting,} iff there exists $g^{\prime }\in \mathrm{Hom}%
_{S}(C,B)$ such that $g\circ g^{\prime }=id_{C}.$

\emph{splitting,} iff it is left splitting and right splitting.
\end{punto}

In 2009, a result to similar Theorem \ref{semi-simple-rings} was proved for
\emph{subtractive} semirings. We add a new characterization using $S$-$e$%
-projective semimodules.

\begin{thm}
\label{id-ss-k-e}If the semiring $S$ is left subtractive, then the following
assertions are equivalent:

\begin{enumerate}
\item Every left $S$-semimodule is $S$-$e$-projective;

\item Every left $S$-semimodule is $S$-$k$-projective;

\item Every short exact sequence $0\rightarrow L\rightarrow S\rightarrow
N\rightarrow 0$ of left $S$-semimodules is right splitting;

\item Every left ideal of $S$ is a direct summand;

\item $S$ is left ideal-semisimple.
\end{enumerate}
\end{thm}

\begin{Beweis}
The equivalences: (2) $\Leftrightarrow $ (4) $\Leftrightarrow $ (5) follow
from \cite[4.4]{KNT2009}.

The equivalences (1) $\Leftrightarrow $ (3) is \cite[Proposition 3.14.]{AN-1}%
.

(1) $\Rightarrow $ (2) follows from the fact that every $S$-$e$-projective
left $S$-semimodule is $S$-$k$-projective.

(4) $\Rightarrow $ (1) This is Lemma \ref{sumproj} applied to $M=$ $%
_{S}S.\blacksquare $
\end{Beweis}

For an arbitrary semiring, having every semimodule projective or injective
or $e$-injective forces the ground semiring to be a semisimple \emph{ring}.
The following observation is a combination of \cite[Theorem 3.1]{Ili2010}
and \cite[5.3]{AIKN2018}:

\begin{thm}
\label{proj-e-inj-ring}The the following assertions are equivalent for any
semiring $S:$

\begin{enumerate}
\item Every left (right) $S$-semimodule is projective;

\item Every left (right) $S$-semimodule is injective;

\item Every left (right) semimodule is $e$-injective;

\item $S$ is a left (right)\ semisimple \textbf{ring}.
\end{enumerate}
\end{thm}

%In 2016, Il'in proved that an injective envelope is not necessary exists.
%Furthermore, he proved the following.

%\begin{thm}
%\cite[3.3]{Ili2016} Every simple $S$-semimodule has injective envelope if
%and only if $S$ is a semizeroic semiring.
%\end{thm}

Our next goal is to find a relationship between the left
ideal-semisimplicity of $S$ and having all left $S$-semimodules $S$-$e$%
-projective.

\begin{defn}
Let $M$ be a left $S$-semimodule. A subsemimodule $N\leq _{S}M$ is a \textbf{%
maximal summand} of $M,$ iff $N\leq _{S}^{\oplus }M$ a direct summand of $M$
such that $N\neq M$ and for every direct summand $L\leq _{S}^{\oplus }M$
with $N\subseteq L\subseteq M,$ we have $N=L$ or $L=M.$ A direct summand $%
N\leq _{S}^{\oplus }M$ is called an \textbf{irreducible summand,} iff $\{0\}$
is a maximal direct summand of $N$.
\end{defn}

\begin{lem}
\label{dcc-acc}\emph{(\cite[Lemma 2.12.]{AN-4})}\ A left $S$-semimodule $M$
satisfies the ACC on direct summands if and only if $M$ satisfies the DCC on
direct summands.
\end{lem}

\begin{lem}
\label{lemint}\emph{(\cite[Lemma 2.3.]{AN-4}) }Let $M$ be an $S$-semimodule
and $N$ a \emph{subtractive }$S$-subsemimodules of $M.$ If $M=L\oplus K$ for
some $L\leq _{S}N$ and $K\leq _{S}M,$ then
\begin{equation*}
N=L\oplus (K\cap N).
\end{equation*}
\end{lem}

\begin{thm}
\label{thmirsum}If $_{S}S$ satisfies the ascending chain condition on direct
summands, then $S=S_{1}\oplus \cdots \oplus S_{n}$, where $S_{i}$ is an
irreducible summand for every $i\in \{1,\cdots ,n\}$.
\end{thm}

\begin{Beweis}
By our assumptions and Lemma \ref{dcc-acc}, $S$ satisfies also the
descending chain condition on direct summands. If $S$ has no non-trivial
direct summand, then $0$ is the maximal summand of $S$, thus $S$ is an
irreducible summand. If not, let $D_{0}$ be a non-trivial direct summand of $%
S$. Then
\begin{equation*}
\mathcal{D}_{1}:=\{D\supsetneqq D_{0}|\text{ }D\text{ is a direct summand of
}S\}
\end{equation*}%
is non-empty as $S\in \mathcal{D}_{1}$. Suppose that there exists $%
(D_{\lambda })_{\Lambda }$ a non-terminating descending chain in $\mathcal{D}%
_{1}$. Then there exists $\lambda _{i}\in \Lambda $, $i=0,1,2,\cdots $ such
that $D_{\lambda _{0}}\supsetneqq D_{\lambda _{1}}\supsetneqq \cdots $, is a
non-terminating strictly descending chain in $\mathcal{D}_{1}$,
contradiction by the DCC on direct summands of $_{S}S$. Thus, the descending
chain $(D_{\lambda })_{\Lambda }$ terminates and has a lower bound.

Since every descending chain in $\mathcal{D}_{1}$ has a lower bound, it
follows by Zorn's Lemma, that $\mathcal{D}_{1}$ has a minimal element, say $%
D_{1}$. Since there is no direct summand between $D_{0}$ and $D_{1}$, we see
that $D_{0}$ is a maximal summand of $D_{1}$.

The set
\begin{equation*}
\mathcal{D}_{-1}:=\{D\subsetneqq D_{0}|\text{ }D\text{ is a direct summand
of }S\}
\end{equation*}%
is non-empty as $0\in \mathcal{D}_{-1}$. Suppose that there exists $%
(D_{\lambda })_{\Lambda }$ a non-terminating ascending chain in $\mathcal{D}%
_{-1}$. Then there exist $\lambda _{i}\in \Lambda $, $i=0,1,\cdots $ such
that $D_{\lambda _{0}}\subsetneqq D_{\lambda _{1}}\subsetneqq \cdots $, is a
non-terminating ascending chain on $\mathcal{D}_{-1}$, contradiction by the
ACC on direct summands of $_{S}S$. Thus the ascending chain $(D_{\lambda
})_{\Lambda }$ terminates and has an upper bound.

Since every ascending chain on $\mathcal{D}_{-1}$ has an upper bound, it
follows by Zorn's Lemma, that $\mathcal{D}_{-1}$ has a maximal element say $%
D_{-1}$. Since there is no direct summand between $D_{-1}$ and $D_{0}$, we
see that $D_{-1}$ is a maximal summand of $D_{0}$. We proved that every
non-trivial direct summand is a maximal summand of a direct summand and has
a maximal summand.

Now, let $D_{0}$ be a non-trivial direct summand of $S.$ Then there exists $%
D_{1},$ a direct summand of $_{S}S,$ such that $D_{0}$ is a maximal summand
of $D_{1}$. If $D_{1}$ is non-trivial, then there exists $D_{2},$ a direct
summand of $S,$ such that $D_{1}$ is a maximal summand of $D_{2}.$ Repeating
this process over and over, we obtain an ascending chain
\begin{equation*}
D_{0}\subsetneqq D_{1}\subsetneqq D_{2}\subsetneqq \cdots
\end{equation*}%
of direct summands of $_{S}S,,$ which should terminate. Thus, there exists $%
n\in \mathbb{N}$ such that
\begin{equation*}
D_{0}\subsetneqq D_{1}\subsetneqq D_{2}\subsetneqq \cdots \subsetneqq D_{n}=S
\end{equation*}%
and $D_{i}$ is maximal summand of $D_{i+1}$ for $i=0,1,...,n-1$. Since $%
D_{0} $ is a non-trivial direct summand of $S$, $D_{0}$ has maximal summand $%
D_{-1} $. If $D_{-1}$ is non-trivial, then $D_{-1}$ has maximal summand $%
D_{-2}$. By repeating this process, we obtain a descending chain
\begin{equation*}
D_{0}\supsetneqq D_{-1}\supsetneqq D_{-2}\supsetneqq \cdots
\end{equation*}%
of direct summands of $_{S}S$, which should terminate. Thus, there exists $%
m\in \mathbb{N}$ such that
\begin{equation*}
D_{0}\supsetneqq D_{-1}\supsetneqq D_{-2}\supsetneqq \cdots \supsetneqq
D_{-m}=0
\end{equation*}%
and $D_{-i}$ is maximal summand of $D_{-i+1}$ for $i=1,2,...,m$. Hence
\begin{equation*}
0=D_{-m}\subsetneqq D_{-m+1}\subsetneqq \cdots \subsetneqq D_{-1}\subsetneqq
D_{0}\subsetneqq D_{1}\subsetneqq D_{2}\subsetneqq \cdots \subsetneqq D_{n}=S
\end{equation*}%
is an ascending chain of direct summands of $S$ such that $D_{i}$ is a
maximal summand of $D_{i+1}$ for $i=-m,-m+1,...,0,1,...,n-1$.

For $i=-m,-m+1,\cdots ,0,1,\cdots ,n-1$, write $S=D_{i}\oplus L_{i}$. Since $%
D_{i}\subsetneqq D_{i+1}$, we have%
\begin{equation*}
D_{i+1}\overset{\text{Lemma \ref{lemint}}}{{=}}{D}_{i}\oplus (D_{i+1}\cap
L_{i}),
\end{equation*}%
with $D_{i+1}\cap L_{i}\neq 0$. Consider $K_{i+1}:=D_{i+1}\cap L_{i}.$ Then
\begin{equation*}
S=D_{n}=K_{-m+1}\oplus K_{-m+2}\oplus \cdots \oplus K_{n}.
\end{equation*}%
Suppose that there exists $i\in \{-m+1,-m+2,\cdots ,n\}$ such that $K_{i}$
is a reducible summand. In this case, there exists a direct summand $K$ of $%
K_{i}$ such that $0\neq K\subsetneqq K_{i}$. Write $K_{i}:=K\oplus L$. Then
\begin{equation*}
S=D_{i}\oplus L_{i}=D_{i-1}\oplus K_{i}\oplus L_{i}=D_{i-1}\oplus K\oplus
L\oplus L_{i},
\end{equation*}%
thus $D_{i-1}\oplus K$ is a direct summand of $S$ such that
\begin{equation*}
D_{i-1}\subsetneqq D_{i-1}\oplus K\subsetneqq D_{i},
\end{equation*}%
contradiction to the maximality of $D_{i-1}$ as summand of $D_{i}$.$%
\blacksquare $
\end{Beweis}

\begin{rem}
\label{rem1}If $S$ is a semiring with $S=\bigoplus\limits_{i\in I}N_{i},$
where $N_{i}$ is a non-zero left ideal of $S$ for every $i\in I$, then $I$
is finite. To see this, suppose that $I$ is infinite. Since $1\in
S=\bigoplus\limits_{i\in I}N_{i}$ we have $1_{S}=\sum%
\limits_{j=1}^{k}n_{i_{j}}$ for some $k\in \mathbb{N},i_{j}\in I$ and $%
n_{i_{j}}\in N_{i_{j}}$. Let $i\in I\backslash \{i_{1},\cdots ,i_{k}\}$ and $%
n_{i}\in N_{i}\backslash \{0\}.$ Then $n_{i}=n_{i}\cdot 1_{S}=n_{i}\cdot
\sum\limits_{j=1}^{k}n_{i_{j}}=\sum\limits_{j=1}^{k}n_{i}n_{i_{j}}$,
contradicting the uniqueness of the representation of $n_{i}$ in the direct
sum.$\blacksquare $
\end{rem}

\begin{prop}
\label{projaccsum}Let $S$ be a semiring such that $S/I$ is $S$-$k$%
-projective for every \emph{subtractive} ideal $I$ of $S.$

\begin{enumerate}
\item $_{S}S$ satisfies the ACC on direct summands.

\item $S=S_{1}\oplus S_{2}\oplus \cdots \oplus S_{n}$, where $S_{i}$ is an
irreducible summand for every $i\in \{1,\cdots ,n\}.$ If moreover, $S_{i}$
is ideal-simple (resp., congruence-simple) for every $i\in\{1,2,...,k\}$,
then $S$ is ideal-semisimple (resp., congruence-semisimple).
\end{enumerate}
\end{prop}

\begin{Beweis}
Assume that $S/I$ is $S$-$k$-projective for every subtractive ideal $I$ of $%
S.$

\begin{enumerate}
\item Suppose, without loss of generality, that there is a strictly
ascending chain of direct summands of $_{S}S$:%
\begin{equation*}
N_{1}\subsetneqq N_{2}\subsetneqq \cdots \subsetneqq N_{i}\subsetneqq
N_{i+1}\subsetneqq \cdots
\end{equation*}%
where, for each $i\in \mathbb{N}$ we have $S=N_{i}\oplus L_{i}$ for some
left ideal $L_{i}\leq _{S}^{\oplus }S.$ Since, $N_{i}\subsetneqq N_{i+1},$we
have $N_{i+1}\overset{\text{Lemma \ref{lemint}}}{{=}}N_{i}\oplus
(N_{i+1}\cap L_{i})$ with $N_{i+1}\cap L_{i}\neq 0$ for each $i\in \mathbb{N}
$. Setting $K_{1}:=N_{1}$ and $K_{i+1}:=N_{i+1}\cap L_{i}$ , for $i\geq 1,$
we have $N_{i}=K_{1}\oplus ...\oplus K_{i}$ for every $i\geq 2$. Thus%
\begin{equation*}
K:=\bigoplus\limits_{i\in \mathbb{N}}K_{i}=\bigcup\limits_{i\in \mathbb{N}%
}N_{i}
\end{equation*}%
is a subtractive left ideal of $S$ as can be easily shown. So, we have an
exact sequence of left $S$-semimodules%
\begin{equation}
0\rightarrow K\overset{\iota }{\longrightarrow }S\overset{\pi }{%
\longrightarrow }S/K\longrightarrow 0.  \label{S/K}
\end{equation}%
Since $S/K$ is $S$-$k$-projective, there exists an $S$-linear map $\varphi
:S/K\rightarrow S$ such that $\pi \circ \varphi =id_{S/K}$. For every $s\in
S $, we have $\pi (\varphi ([1]))=(\pi \circ \varphi )([1])=[1].$ Sine $\pi $
is $k$-normal, there exist $k,k^{\prime }\in K$ such that $1+k=\varphi
([1])+k^{\prime }.$ Write $k=k_{1}+\cdots +k_{j}$ and $k^{\prime
}=k_{1}^{\prime }+\cdots +k_{l}^{\prime },$ where $k_{i},k_{i}^{\prime }\in
K_{i}$ for every $i,$ and let $m:=max\{j,l\}$. Then $k=k_{0}+k_{1}+...+k_{m}$
and $k^{\prime }=k_{1}^{\prime }+...+k_{m}^{\prime }$ for some $%
k_{i},k_{i}^{\prime }\in K_{i}$. Recall that for every $i\in \mathbb{N}$ we
have%
\begin{equation*}
S=N_{i}\oplus L_{i}=(K_{1}\oplus ...\oplus K_{i-1})\oplus K_{i}\oplus L_{i}.
\end{equation*}%
For every $i\in \mathbb{N}$, let $\pi _{i}:S\rightarrow K_{i}$ be the
canonical projection on $K_{i}$ and $e_{i}:=\pi _{i}(1).$ Then, $%
e_{i}=e_{i}1 $ implies $\pi _{j}(e_{i})=\pi _{j}(e_{i}1)=e_{i}\pi
_{j}(1)=e_{i}e_{j}$ and so $e_{i}e_{j}=0$ for every $i\neq j$ and $%
e_{i}e_{i}=e_{i}$. Since $k,k^{\prime }\in N_{m}$, we have $\pi
_{m+1}(k)=0=\pi _{m+1}(k^{\prime })$. Thus
\begin{equation*}
e_{m+1}=\pi _{m+1}(1+k)=\pi _{m+1}(\varphi ([1])+k^{\prime })=\pi
_{m+1}(\varphi ([1])).
\end{equation*}%
Since $S=N_{m+1}\oplus L_{m+1}=K_{1}\oplus ...\oplus K_{m}\oplus
K_{m+1}\oplus L_{m+1}$, we have
\begin{equation*}
1=e_{1}+...+e_{m}+e_{m+1}+l_{m+1}=(e_{1}+...+e_{m}+l_{m+1})+e_{m+1}
\end{equation*}%
for some $l_{m+1}\in L_{m+1}$, whence
\begin{equation*}
\pi (1)=\pi (e_{0}+e_{1}+...+e_{m}+l_{m+1}),\text{ i.e. }%
[1]=[e_{1}+...+e_{m}+l_{m+1}].
\end{equation*}

Notice that%
\begin{equation*}
\begin{tabular}{lll}
$\varphi ([e_{1}+...+e_{m}+l_{m+1}])$ & $=$ & $\varphi (\pi
(e_{0}+e_{1}+...+e_{m}+l_{m+1}))$ \\
& $=$ & $\varphi (\pi ((e_{1}+...+e_{m}+l_{m+1})1))$ \\
& $=$ & $\varphi ((e_{1}+...+e_{m}+l_{m+1})\pi (1))$ \\
& $=$ & $(e_{1}+...+e_{m}+l_{m+1})\varphi (\pi (1))$ \\
& $=$ & $(e_{0}+e_{1}+...+e_{m}+l_{m+1})\varphi ([1])$ \\
& $=$ & $\pi _{m+1}(\varphi ([e_{0}+e_{1}+...+e_{m}+l_{m+1}]))$ \\
& $=$ & $\pi _{m+1}((e_{0}+e_{1}+...+e_{m}+l_{m+1})\varphi ([1]))$ \\
& $=$ & $(e_{0}+e_{1}+...+e_{m}+l_{m+1})\pi _{m+1}(\varphi ([1]))$ \\
& $=$ & $(e_{0}+e_{1}+...+e_{m}+l_{m+1})e_{m+1}$ \\
& $=$ & $l_{m+1}e_{m+1}$ \\
& $=$ & $l_{m+1}\pi _{m+1}(1)$ \\
& $=$ & $\pi _{m+1}(l_{m+1})$ \\
& $=$ & $0.$%
\end{tabular}%
\end{equation*}%
It follows that $[1]=[e_{0}+e_{1}+...+e_{m}+l_{m+1}]$ while $\varphi
([e_{0}+e_{1}+...+e_{m}+l_{m+1}])\neq \varphi ([1])$, a contradiction.
Hence, every ascending chain of direct summands of $S$ terminates, i.e. $S$
satisfies the ACC on direct summands.

\item By (1), the assumptions of Theorem \ref{thmirsum} are satisfied,
whence
\begin{equation*}
S=S_{1}\oplus \cdots \oplus S_{n}
\end{equation*}%
where $S_{i}$ is an irreducible summand for every $i\in \{1,\cdots ,n\}$. If
moreover, $S_{i}$ is ideal-simple (resp., congruence-simple) for every $i\in
\{1,\cdots ,n\}$, then $S$ is the direct sum of ideal-simple (resp.
congruence-simple) left ideals, whence ideal-semisimple (resp.
congruence-semisimple).$\blacksquare $
\end{enumerate}
\end{Beweis}

The following result is a combination of Lemma \ref{sumproj} and Proposition %
\ref{projaccsum}.

\begin{cor}
\label{sumidsim}If $S$ is a semiring such that every subtractive left ideal
is a direct summand, then $S=S_{1}\oplus \cdots \oplus S_{n}$, where $S_{i}$
is an irreducible summand for every $i\in \{1,\cdots ,n\}.$ If moreover, $%
S_{i}$ is ideal-simple (resp., congruence-simple) for every $i\in \{1,\cdots
,n\}$, then $S$ is ideal-semisimple (resp., congruence-semisimple).
\end{cor}

\begin{prop}
\label{proj-impl}For any semiring $S,$ each of the following conditions
implies its successor:

\begin{enumerate}
\item Every subtractive ideal of $S$ is a direct summand.

\item Every $S$-semimodule is $S$-$e$-projective.

\item Every $S$-semimodule is $S$-$k$-projective.

\item $S/I$ is $S$-$k$-projective for every subtractive ideal $I$ of $S.$

\item Every short exact sequence $0\longrightarrow I\longrightarrow
S\longrightarrow N\longrightarrow 0$ in $_{S}\mathbf{SM}$ right splits.

\item $_{S}S$ satisfies ACC on direct summands.

\item $_{S}S$ satisfies DCC on direct summands.

\item $S=S_{1}\oplus \cdots \oplus S_{n}$, where every $S_{i}$ is an
irreducible summand.
\end{enumerate}
\end{prop}

\begin{Beweis}
(1) $\Rightarrow $ (2) This follows from Lemma \ref{sumproj} applied to $M=$
$_{S}S$. Let $M$ be an irreducible summand of $_{S}S,$ i.e. $\{{0\}}$ is the
only maximal direct summand of $_{S}M$. By our assumption, $M\simeq M/{0}$
is ideal-simple.

(2) $\Rightarrow $ (3) $\Rightarrow $ (4) Follow directly from the
definitions.

(4) $\Longleftrightarrow $ (5 ) Follows from \cite[Proposition 3.14]{AN-1}
and Lemma \ref{exact}.

(4) $\Rightarrow $ (6) Follows from Proposition \ref{projaccsum}.

(6) $\Leftrightarrow $ (7) Follows from Lemma \ref{dcc-acc}.

(6) $\Rightarrow $ (8) Follows by Theorem \ref{thmirsum}.$\blacksquare $
\end{Beweis}

\begin{thm}
\label{left-spit-k-Noeth}\emph{(\cite[Theorem 2.21]{AN-4})} If $S$ is a
semiring such that every short exact sequence of left $S$-semimodules $%
0\rightarrow L\rightarrow S\rightarrow N\rightarrow 0$ is left splitting,
then $S$ is a left $k$-Noetherian.
\end{thm}

\begin{prop}
\label{sum-einj}For any semiring $S,$ each of the following conditions
implies its successor:

\begin{enumerate}
\item Every subtractive left ideal of $S$ is a direct summand.

\item Every left $S$-semimodule is $S$-$e$-injective.

\item Every $S$-semimodule is $S$-$i$-injective.

\item Every subtractive ideal of $S$ is $S$-$i$-injective.

\item Every short exact sequence $0\rightarrow L\rightarrow S\rightarrow
N\rightarrow 0$ in $_{S}\mathbf{SM}$ is left splitting.

\item $S$ is $k$-Noetherian.

\item $S$ satisfies the ACC on direct summands.

\item $S$ satisfies the DCC on direct summands.

\item $S=S_{1}\oplus S_{2}\oplus S_{3}\oplus ...\oplus S_{n}$, where every $%
S_{i}$ is an irreducible summand.
\end{enumerate}
\end{prop}

\begin{Beweis}
(1) $\Rightarrow $ (2) Let $J$ be a left $S$-semimodule and let $%
f:M\rightarrow S$ a normal monomorphism, i.e. $M\leq S$ is a subtractive
left ideal and $f$ is the canonical embedding. Let $g:M\rightarrow J$ be an $%
S$-linear map. By the assumption, $S=M\oplus N$ for some left ideal $N$ of $%
S.$ Let $\pi :S\rightarrow M$ be the projection on $M$ (i.e., $\pi \circ
f=id_{M}$). Then $g\circ \pi :S\rightarrow J$ is an $S$-linear map
satisfying $(g\circ \pi )\circ f=g.$

Let $h:M\rightarrow J$ be another $S$-linear map satisfying $h\circ f=g.$
Write $1_{S}=e_{M}+e_{N}$, where $e_{M}\in M$ are $e_{N}\in N$ are uniquely
determined, and let $j_{0}:=h(1_{S}).$ For every $m\in M,$ we have $%
m=m1_{S}=m(e_{M}+e_{N})=me_{M}+me_{N}$, whence $me_{M}=m$ and $me_{N}=0$ as
the sum $M+N$ is direct. Similarly, $ne_{M}=0$ and $ne_{N}=n$ for every $%
n\in N$. Define
\begin{equation*}
h_{1}:S\rightarrow J,\text{ }s\mapsto se_{N}j_{0}.
\end{equation*}%
Then $(h_{1}\circ f)(m)=h_{1}(m)=me_{N}j_{0}=0$ for every $m\in M$.
Moreover, we have%
\begin{equation*}
\begin{array}{ccccc}
(g\circ \pi +h_{1})(s) & = & (g\circ \pi )(s)+h_{1}(s) & = & (g\circ \pi
)(se_{M}+se_{N})+h_{1}(s) \\
& = & g(se_{M})+se_{N}j_{0} & = & (h\circ f)(se_{M})+se_{N}j_{0} \\
& = & h(se_{M})+se_{N}j_{0} & = & se_{M}j_{0}+se_{N}j_{0} \\
& = & s(e_{M}+e_{N})j_{0} & = & sj_{0} \\
& = & h(s) & = & (h+0)(s).%
\end{array}%
\end{equation*}%
Hence $J$ is $S$-$e$-injective.

The implications (2) $\Rightarrow $ (3) $\Rightarrow $ (4) $\Rightarrow $
(5) $\&\ $(6) $\Rightarrow $ (7) follow from the definitions.

(5) $\Rightarrow $ (6) Follows from Theorem \ref{left-spit-k-Noeth}.

(7) $\Longleftrightarrow $ (8) Follows from Lemma \ref{dcc-acc}.

(7) $\Rightarrow $ (9) Follows from Theorem \ref{thmirsum}.$\blacksquare $
\end{Beweis}

\section{Commutative semisimple Semirings}

The converse of Corollary \ref{sumidsim} is satisfied when the semiring $S$
is commutative. To achieve this, we first prove the following technical
result.

\begin{lem}
\label{comsum}Let $S$ be a \emph{commutative} ideal-semisimple
(congruence-semisimple) semiring and write $S=S_{1}\oplus S_{2}\oplus
...\oplus S_{k},$ where $S_{i}$ is an ideal-simple ideal of $S$ for every $%
i\in \{1,,\cdots ,k\}$. Then every subtractive ideal $I$ of $S$ is a direct
summand, and moreover $I=\bigoplus\limits_{a\in A}S_{a}$ for some $%
A\subseteq \{1,\cdots ,k\}$.
\end{lem}

\begin{Beweis}
Let $I$ be a subtractive ideal of $S$ and
\begin{equation*}
A=\{a\in \{1,\cdots ,k\}|\text{ }I\cap S_{a}\neq \{0\}\}.
\end{equation*}%
Let $B:=\{1,\cdots ,k\}\backslash A$ and write $S_{A}:=\bigoplus\limits_{a%
\in A}S_{a}$ and $S_{B}:=\bigoplus\limits_{b\in B}S_{b}$ . For every $a\in A$%
, the ideal $S_{a}$ is a (subtractive) ideal of $A,$ thus $I\cap S_{a}$ is a
(subtractive) ideal. Since $0\neq I\cap S_{a}\subseteq S_{a}$ and $I\cap
S_{a}$ is a (subtractive) left ideal, $I\cap S_{a}=S_{a}$. Thus $%
S_{A}\subseteq I$, and it follows that $I\overset{\text{Lemma \ref{lemint}}}{%
{=}}S_{A}\oplus (S_{B}\cap I)$.

\textbf{Claim: }$I\cap S_{B}=0$.

Let $1=e_{1}+...+e_{k}$ for some $e_{i}\in S_{i}$. For every $s_{i}\in S_{i}$%
,
\begin{equation*}
s_{i}=s_{i}1=s_{i}(e_{1}+e_{2}+...+e_{k})=s_{i}e_{1}+s_{i}e_{2}+...+s_{i}e_{k}.
\end{equation*}%
Since $s_{i}e_{j}\in S_{j}$ for every $j\in \{1,,\cdots ,k\}$, it follows by
the directness of the sum that $s_{i}e_{i}=s_{i}$ and $s_{i}e_{j}=0$ for
every $i\neq j$. Therefore $e_{i}s_{i}=s_{i}$ and $e_{j}s_{i}=0$ for every $%
i\neq j$. Let $x\in I\cap S_{B}$, whence $x=\sum\limits_{b\in B}x_{b}$ where
$x_{b}\in S_{b}$ for each $b\in B$. For every $\tilde{b}\in B$, we have $x_{%
\tilde{b}}=\sum\limits_{b\in B}e_{\tilde{b}}x_{b}=e_{\tilde{b}}x\in I$ as $I$
is an ideal. Thus $x_{b}=0$ for every $b\in B$ and $x=0$.$\blacksquare $
\end{Beweis}

The following technical conditions shall be needed in the sequel.

\begin{punto}
\label{C1-C2}Let $N$ be a left $S$-semimodule. Consider the conditions:

$\mathbf{C1:\ }$Every subtractive $S$-subsemimodule $M\leq_{S} N$ is a
direct summand.

$\mathbf{C2:\ }$For every subtractive $S$-subsemimodule $M\leq _{S}N$ and
every maximal subtractive $S$-subsemimodule $L\leq_{S} M$, the left $S$%
-semimodule $M/L$ is left ideal-simple.

$\mathbf{C2}^{\prime }\mathbf{:}$ For every subtractive subsemimodule $M\leq
_{S} N$ and every maximal subtractive $S$-subsemimodule $L \leq_{S} M$, the
left $S$-semimodule $M/L$ is congruence-simple.
\end{punto}

\begin{rem}
\label{indp}The conditions $\mathbf{C1}$ and $\mathbf{C2}$ (and $\mathbf{C2}%
^{\prime }$) are independent:

\begin{enumerate}
\item $B(3,2)$ satisfies $\mathbf{C1}$ but neither $\mathbf{C2}$ nor $%
\mathbf{C2}^{\prime }.$

\item $B(3,1)$ satisfies $\mathbf{C2}$ but not $\mathbf{C1}.$

\item $\mathbb{B}^{\mathbb{N}}$ satisfies $\mathbf{C2^{\prime }}$ but not $%
\mathbf{C1}.$

\item $\mathbb{R}^{+}$ satisfies $\mathbf{C2}$ but not $\mathbf{C2^{\prime }.%
}$ By Example \ref{id-not-cong}, $\mathbb{R}^{+}$ is ideal-simple but not
congruence-simple. Since $\mathbb{R}^{+}$ is ideal-simple, it has no proper
non-trivial ideals, $\{0\}$ is the maximal subtractive subsemimodule of $%
\mathbb{R}^{+},$ and $\mathbb{R}^{+}/\{0\}\simeq \mathbb{R}^{+}$ is
ideal-simple. Hence $\mathbb{R}^{+}$ satisfies $\mathbf{C2}$. However, $%
\mathbb{R}^{+}/\{0\}\simeq \mathbb{R}^{+}$ is not congruence-simple, thus $%
\mathbb{R}^{+}$ does not satisfy $\mathbf{C2^{\prime }.}$

\item Let $(M,+,0)$ be a finite lattice which is not distributive. $\mathbf{E%
}_{M},$ the endomorphism semiring of $M,$ satisfies $\mathbf{C2^{\prime }}$
but not $\mathbf{C2.}$ By Example \ref{cong-not-id}, $\mathbf{E}_{M}$ is
left congruence-simple but not left ideal-simple. Since $\mathbf{E}_{M}$ is
left congruence-simple, it has no non-trivial subtractive left ideals, $%
\{0\} $ is the maximal subtractive ideal of $\mathbf{E}_{M}$ and $\mathbf{E}%
_{M}/\{0\}=\mathbf{E}_{M}$ is left congruence-simple. Hence, $\mathbf{E}_{M}$
satisfies $\mathbf{C2^{\prime }}$. However, $\mathbf{E}_{M}/\{0\}=\mathbf{E}%
_{M}$ is not ideal-simple, thus $\mathbf{E}_{M}$ does not satisfy $\mathbf{C2%
}$.
\end{enumerate}
\end{rem}

The following result extends the characterizations of ideal-semisimple
semirings in Theorem \ref{id-ss-k-e} to commutative \emph{not necessarily
subtractive} semirings:

\begin{thm}
\label{id-SS-C1}The following assertions are equivalent for a \emph{%
commutative} semiring $S$:

\begin{enumerate}
\item Every subtractive ideal of $S$ is a direct summand and $S$ satisfies $%
\mathbf{C2};$

\item Every $S$-semimodule is $S$-$e$-projective and $S$ satisfies $\mathbf{%
C2};$

\item Every $S$-semimodule is $S$-$k$-projective and $S$ satisfies $\mathbf{%
C2};$

\item $S/I$ is $S$-$k$-projective for every subtractive ideal $I$ of $S,$
and $S$ satisfies $\mathbf{C2};$

\item Every short exact sequence $0\longrightarrow I\longrightarrow
S\longrightarrow N\longrightarrow 0$ in $_{S}\mathbf{SM}$ right splits and $%
S $ satisfies $\mathbf{C2};$

\item $_{S}S$ satisfies the ACC on direct summands and $\mathbf{C2};$

\item $_{S}S$ satisfies the DCC on direct summands and $\mathbf{C2};$

\item $S=S_{1}\oplus S_{2}\oplus S_{3}\oplus ...\oplus S_{n}$, where every $%
S_{i}$ is an irreducible summand, and $S$ satisfies $\mathbf{C2};$

\item $S$ is ideal-semisimple.
\end{enumerate}
\end{thm}

\begin{Beweis}
By Proposition \ref{proj-impl}, we only need to prove (8) $\Rightarrow $ (9)
and (9) $\Rightarrow $ (1).

(8) $\Rightarrow $ (9) Notice that assuming $\mathbf{C2}$ guarantees that $%
S_{i}$ is ideal-simple for $i=1,\cdots ,n.$ Whence, $S$ is ideal-semisimple.

(9) $\Rightarrow $ (1) Assume that $S$ is ideal-semisimple and write $%
S=S_{1}\oplus \cdots \oplus S_{k}$ for some $k\in \mathbb{N}$ with $S_{i}$
an ideal-simple ideal for $i=1,\cdots ,k$. Let $I$ be a subtractive ideal of
$S$. Since $S$ is commutative, it follows by Lemma \ref{comsum} that $%
I=\bigoplus\limits_{a\in A}S_{a}$ for some $A\subseteq \{1,\cdots ,k\}$,
whence $S\overset{\text{Lemma \ref{lemint}}}{{=}}I\oplus
\bigoplus\limits_{b\notin A}S_{b}$. Hence, $I$ is a direct summand of $S$.

\textbf{Claim:} $_{S} S$ satisfies ${\mathbf{C2}}$.

Let $M$ be a subtractive ideal of $S$ and $L$ a maximal subtractive subideal
of $M$. Then $M=S_{A}=\bigoplus\limits_{a\in A}S_{a}$ and $%
L=S_{C}=\bigoplus\limits_{c\in C}S_{c}$ for some $C\subsetneqq A\subseteq
\{1,\cdots ,k\}$. Notice that $C\subseteq A$ since $L\subseteq M$. Moreover,
$|A\backslash C|=1$ since $|A\backslash C|=0$ implies $L=M$ and $\left\vert
A\backslash C\right\vert \geq 2$ implies $L\subsetneqq S_{C\cup
\{y\}}\subsetneqq M$ for some $y\in A\backslash C$ with $S_{C\cup \{y\}}$ a
subtractive ideal of $S$, contradiction to the maximality of $L$. Write $%
A\backslash C=\{x\}$ and $B=\{1,2,...,k\}\backslash A$. Then $S=S_{A}\oplus
S_{B}=S_{C}\oplus S_{x}\oplus S_{B}$ where $S_{B}=\bigoplus\limits_{b\in
B}S_{b}$.

Let $I$ be an ideal of $S$ such that $L\subsetneqq I\subseteq M.$ Then there
exists $i\in I\backslash N$. Since $i\in M$, $i=t_{C}+t_{x}$ for some $%
t_{C}\in S_{C},t_{x}\in S_{x}$. Notice that $t_{x}\neq 0$; otherwise, $%
i=t_{C}\in N$. Moreover, $0\neq
t_{x}=e_{x}t_{x}=e_{x}(t_{C}+t_{x})=e_{x}i\in I$. Thus $I\cap S_{x}\neq 0$,
whence $I\cap S_{x}=S_{x}$ as $S_{x}$ is ideal-simple. Since $S_{C}\subseteq
I$ and $S_{x}\subseteq I$, we have $M=S_{C}+S_{x}\subseteq I$.$\blacksquare $
\end{Beweis}

The following result is the ``congruence-semisimple" version of Theorem \ref%
{id-SS-C1}.

\begin{thm}
\label{cong-C2}The following assertions are equivalent for a {commutative}
semiring $S$:

\begin{enumerate}
\item Every subtractive ideal of $S$ is a direct summand and $S$ satisfies $%
\mathbf{C2}^{\prime };$

\item Every $S$-semimodule is $S$-$e$-projective and $S$ satisfies $\mathbf{%
C2}^{\prime };$

\item Every $S$-semimodule is $S$-$k$-projective and $S$ satisfies $\mathbf{%
C2}^{\prime };$

\item $S/I$ is $S$-$k$-projective for every subtractive ideal $I$ of $S$ and
$S$ satisfies $\mathbf{C2}^{\prime };$

\item Every short exact sequence $0\longrightarrow I\longrightarrow
S\longrightarrow N\longrightarrow 0$ right splits and $S$ satisfies $\mathbf{%
C2}^{\prime };$

\item $_{S}S$ satisfies the ACC on direct summands and $S$ satisfies $%
\mathbf{C2}^{\prime };$

\item $_{S}S$ satisfies the DCC on direct summands and $S$ satisfies $%
\mathbf{C2}^{\prime };$

\item $S=S_{1}\oplus \cdots \oplus S_{n}$, where every $S_{i}$ is an
irreducible summand, and $S$ satisfies $\mathbf{C2}^{\prime };$

\item $S$ is congruence-semisimple.
\end{enumerate}
\end{thm}

\begin{Beweis}
We only need to prove (9) $\Rightarrow $ (1); the proof of the other
implications are similar to the proof of the corresponding ones in Theorem %
\ref{id-SS-C1}.

Assume that $S$ is congruence-semisimple. With the help of Lemma \ref{comsum}%
, it can be shown, as in the proof of Theorem \ref{id-SS-C1}, that every
subtractive ideal of $S$ is a direct summand.

\textbf{Claim:\ }$S$ satisfies $\mathbf{C2}^{\prime }\mathbf{.}$

Let $M,$ $L$ be subtractive ideals of $S$ with $L$ a maximal subtractive $S$%
-subsemimodule of $M$. Then similarly to the proof of Theorem \ref{id-SS-C1}%
, we have $M=S_{A}:=\bigoplus\limits_{a\in A}S_{a}$, $S=S_{A}\oplus S_{B}$
and $L=S_{C}:=\bigoplus\limits_{c\in C}S_{c}$ where $C\cup \{x\}=A$.

Let $\rho $ be a congruence relation on $S$ such that $\equiv
_{L}\subsetneqq \rho \subseteq \equiv _{M}.$ Consider the congruence
relation $\rho ^{\prime }$ on $S_{x}$ defined by%
\begin{equation*}
t_{x}\rho ^{\prime }t_{x}^{\prime }\Leftrightarrow (t_{C}+t_{x}+t_{B})\rho
(t_{C}^{\prime }+t_{x}^{\prime }+t_{B}^{\prime })\text{ for some }%
t_{C},t_{C}^{\prime }\in S_{C},t_{B},t_{B}^{\prime }\in S_{B}.
\end{equation*}

\textbf{Step I:\ } $\rho ^{\prime }=S_{x}^{2}$.

Since $\equiv _{N}\neq \rho $, there exist $s,s^{\prime }\in S$ such that $%
s\not\equiv _{L}s^{\prime }$ and $s\rho s^{\prime }$. Write $%
s=s_{C}+s_{x}+s_{B}$ and $s^{\prime }=s_{C}^{\prime }+s_{x}^{\prime
}+s_{B}^{\prime }$ for some $s_{C},s_{C}^{\prime }\in
S_{C},s_{x},s_{x}^{\prime }\in S_{x},s_{B},s_{B}^{\prime }\in S_{B}$. Since $%
s\equiv _{M}s^{\prime }$, there exists $m,m^{\prime }\in M=S_{A}$ such that $%
m+s=m^{\prime }+s^{\prime }$, that is $(m+s_{C}+s_{x})+s_{B}=(m^{\prime
}+s_{C}^{\prime }+s_{x}^{\prime })+s_{B}^{\prime }$ whence $%
s_{B}=s_{B}^{\prime }$ as the sum $S_{A}+S_{B}$ is direct. Notice that $%
s_{x}\neq s_{x}^{\prime }$; otherwise, $s_{C}^{\prime }+s=s_{C}+s^{\prime }$
where $s_{C},s_{C}^{\prime }\in S_{C}=L$, a contradiction (with $\equiv
_{L}\neq \rho $). Therefore, $s_{x}\rho ^{\prime }s_{x}^{\prime }$ and $%
s_{x}\neq s_{x}^{\prime }$, whence $\rho ^{\prime }=S_{x}^{2}$ as $S_{x}$ is
congruence-simple.

\textbf{Step II:\ } $\rho =\equiv _{M}$.

Let $s,s^{\prime }\in S$ be such that $s\equiv _{M}s^{\prime }$ and write $%
s=s_{C}+s_{x}+s_{B}$, $s^{\prime }=s_{C}^{\prime }+s_{x}^{\prime
}+s_{B}^{\prime }$ for some $s_{C},s_{C}^{\prime }\in
S_{C},s_{x},s_{x}^{\prime }\in S_{x},s_{B},s_{B}^{\prime }\in S_{B}$. Then $%
s_{B}=s_{B}^{\prime }$. Since $\rho ^{\prime }=S_{x}^{2}$, we have $%
s_{x}\rho ^{\prime }s_{x}^{\prime }$, whence $(t_{C}+s_{x}+t_{B})\rho
(t_{C}^{\prime }+s_{x}^{\prime }+t_{B}^{\prime })$ for some $%
t_{C},t_{C}^{\prime }\in S_{C},t_{B},t_{B}^{\prime }\in S_{B}$. Thus $%
e_{x}(t_{C}+s_{x}+t_{B})\rho e_{x}(t_{C}^{\prime }+s_{x}^{\prime
}+t_{B}^{\prime })$, that is $s_{x}\rho s_{x}^{\prime }$. Since $s_{C}\equiv
_{L}s_{C}^{\prime },$ we have $s_{x}\rho s_{x}^{\prime }$, and $%
s_{B}=s_{B}^{\prime }$, $(s_{C}+s_{x}+s_{B})\rho (s_{C}^{\prime
}+s_{x}^{\prime }+s_{B}^{\prime })$, that is $s\rho s^{\prime }$. We
conclude that $\rho =\equiv _{M}$.$\blacksquare $
\end{Beweis}

\begin{definition}
The semiring $S$ is

\textbf{left }(\textbf{right})\textbf{\ }$k$\textbf{-Noetherian}, iff every
ascending chain of \emph{subtractive} left (right) ideals of $S$ terminates;

\textbf{left }(\textbf{right})\textbf{\ }$k$\textbf{-Artinian}, iff every
descending chain of \emph{subtractive} left (right) ideals of $S$ terminates.
\end{definition}

\begin{thm}
\label{dss-sArt}\emph{(\cite[Theorem 2.13]{AN-4})} If every subtractive left
ideal of $S$ is a direct summand, then $S$ is left $k$-Artinian and left $k$%
-Noetherian.
\end{thm}

The following result is a combination of Theorems \ref{id-SS-C1}, \ref%
{cong-C2} and \ref{dss-sArt}:

\begin{cor}
\label{c-ss-N-A}If $S$ is a \emph{commutative} ideal-semisimple
(congruence-semisimple) semiring, then $S$ is $k$-Artinian and $k$%
-Noetherian.
\end{cor}

The following examples show that $\mathbf{C2}$ (resp., $\mathbf{C2}^{\prime
} $) cannot be dropped from the assumptions of Theorem \ref{id-SS-C1}
(resp., Theorem \ref{cong-C2}).

\begin{ex}
\label{exb32}Consider the commutative semiring $B(p+1,p),$ where $p$ is an
odd prime number.

\begin{enumerate}
\item Every subtractive ideal is a direct summand.

\item Every $B(p+1,p)$-semimodule is $B(p+1,p)$-$e$-projective.

\item $B(p+1,p)$ is not left ideal-semisimple.

\item $B(p+1,p)$ is not congruence-semisimple.
\end{enumerate}
\end{ex}

\begin{Beweis}
Notice that the only ideals of $B(p+1,p)$ are $\{0\},$ $B(p+1,p)$, and $%
I=\{0,p\}.$

\begin{enumerate}
\item The only subtractive ideals of $B(p+1,p)$ are $\{0\}$ and $B(p+1,p),$
each of which is a direct summand of $B(p+1,p).$

\item Since (1) is valid, it follows by Lemma \ref{sumproj}, that all $%
B(p+1,p)$-ideals are $B(p+1,p)$-$e$-projective.

\item $B(p+1,p)$ is an irreducible summand which is not ideal-simple (it
contains the ideal $I$). So, $B(p+1,p)$ is not ideal-semisimple. Notice that
$B(p+1,p)$ does not satisfy $\mathbf{C2}$.

\item $B(p+1,p)$ is not an irreducible summand which is not
congruence-simple ($\rho =\{(i,j)|$ $i,j\neq 0\}$ is a non trivial
congruence relation on $B(p+1,p)$). Notice that $B(p+1,p)$ does not satisfy $%
\mathbf{C2}^{\prime }$.$\blacksquare $%
%
%
%
%
%
%
%
%
%
%
%
%
%
%
%
%
%
%
%
%
%
%
%
%
%
%
%
%
%
%
%
%
%
% is subtractive-Artinian where every subtractive left
%ideal of $S$ is a direct summand. But $S$ is not a left ideal-semisimple
%semiring, since it is not left ideal-simple as $I=\{0,2\}$ is a left ideal
%of $S$ with $0\subsetneqq I\subsetneqq S$. Moreover $0$ is a maximal
%subtractive $S$-subsemimodules of $S$ but $S/0\cong S$ is not ideal simple
%as $0\lneqq _{S}I\lneqq _{S}S$.
\end{enumerate}
\end{Beweis}

\begin{ex}
\label{B(3,1)}Consider the semiring $S:=B(3,1).$

\begin{enumerate}
\item $I:=\{0,2\}$ is a subtractive ideal of $B(3,1)$, which is not a direct
summand of $B(3,1)$;

\item $B(3,1)$ is not ideal-semisimple;

\item $B(3,1)$ is not congruence-semisimple.
\end{enumerate}
\end{ex}

\begin{Beweis}
Notice that the only ideals of $S$ are $0,$ $I$ and $S,$ which are
subtractive. Moreover, $I$ is the maximal subtractive subsemimodule of $S$
and is clearly not a direct summand of $S.$ Moreover, $\{0_{S}\}$ is the
maximal subtractive ideal of $I.$ Notice that $I/0\cong \mathbb{B}\cong S/I$
as $S$-semimodules, whence $I/0$ and $S/I$ are ideal-simple. Thus $S$ is an
irreducible summand that is neither ideal-simple ($I$ is a non trivial left
ideal of $S$ ) nor congruence-simple ($\equiv _{I}$ is a non trivial
congruence relation of $S$).
%Recall that $B(3,1)=\{0,1,2\}$ with $2+1=1$. Let $S=B(3,1)$, then $S$ is
%subtractive-Artinian as it has finitely many elements. The only left ideals
%of $S$ are $0,I=\{0,2\}$, and $S$, all of them are subtractive. $S$ is not
%congruence-semisimple as it has only three elements and $S$ is not
%congruence-simple. $I$ is not a direct summand, since whenever a left ideal $%
%J$ satisfies $I+J=S$, $J$ must contains $1$, which implies $J=S$ and the sum
%$I+J$ is not direct. Note that $I/0\cong I$ and $I$ is congruence-simple. $%
%S/I\cong \mathbb{B}$ and $\mathbb{B}$ is congruence-simple.
%$\equiv_0=\Delta_S=\{(0,0),(1,1),(2,2)\}$, $\equiv_I=\{(0,0),(1,1),(2,2),(0,2),(2,0)\}$, and $\equiv_S=S^2$. If $\rho$ is a congruence relation between $\equiv_0$ and $\equiv_I$ with $\rho\neq\equiv_0$, then $0\rho 2$, therefore $\rho=\equiv_I$. If $\rho'$ is a congruence relation between $\equiv_I$ and $\equiv_S$ with $\rho'\neq\equiv_I$, then $1\rho' 2$ or $1\rho' 0$, either case implies $\rho'=\equiv_S$. Hence whenever $M$ is a non-zero subtractive ideal of $S$ and $N$ is a maximal subtractive subsemimodule of $N$, there is no congruence relation between $\equiv_N$ and $\equiv_M$.
\newline
\end{Beweis}

\begin{thm}
\label{ISS-COMM}Let $S$ be a \emph{commutative} semiring which satisfies $%
\mathbf{C2}.$ The following assertions are equivalent:

\begin{enumerate}
\item Every subtractive ideal of $S$ is a direct summand;

\item Every $S$-semimodule is $S$-$e$-injective;

\item Every $S$-semimodule is $S$-$i$-injective;

\item Every subtractive ideal of $S$ is $S$-$i$-injective;

\item Every short exact sequence $0\rightarrow L\rightarrow S\rightarrow
N\rightarrow 0$ in $_{S}\mathbf{SM}$ is left splitting;

\item $S$ is $k$-Noetherian;

\item $S$ satisfies the ACC on direct summands;

\item $S$ satisfies the DCC on direct summands;

\item $S=S_{1}\oplus S_{2}\oplus S_{3}\oplus ...\oplus S_{n}$, where every $%
S_{i}$ is an irreducible summand;

\item $S$ is ideal-semisimple.
\end{enumerate}
\end{thm}

\begin{Beweis}
This is a consequence of Proposition \ref{sum-einj} and the proof of Theorem %
\ref{id-SS-C1}.$\blacksquare $
\end{Beweis}

Combining Theorems \ref{id-SS-C1} and \ref{ISS-COMM}, we obtain the
following characterization of commutative ideal-semisimple semirings:

\begin{thm}
\label{com-idss}The following assertions are equivalent for a commutative
semiring $S$ which satisfies $\mathbf{C2}$:

\begin{enumerate}
\item Every subtractive ideal of $S$ is a direct summand;

\item Every $S$-semimodule is $S$-$e$-projective ($S$-$k$-projective);

\item Every $S$-semimodule is $S$-$e$-injective ($S$-$i$-injective);

\item For every subtractive ideal $I$ of $S$ we have: $S/I$ is $S$-$k$%
-projective ($I$ is $S$-$i$-injective);

\item Every short exact sequence $0\longrightarrow I\longrightarrow
S\longrightarrow N\longrightarrow 0$ in $_{S}\mathbf{SM}$ right splits (left
splits);

\item $S$ is $k$-Noetherian;

\item $_{S}S$ satisfies ACC on the direct summands;

\item $_{S}S$ satisfies DCC on the direct summands;

\item $S=S_{1}\oplus S_{2}\oplus S_{3}\oplus ...\oplus S_{n}$, where every $%
S_{i}$ is an irreducible summand;

\item $S$ is ideal-semisimple.
\end{enumerate}
\end{thm}

The following result is the congruence-semisimple version of Theorem \ref%
{ISS-COMM}.

\begin{thm}
\label{CSS-COMM}Let $S$ be a commutative semiring which satisfies $\mathbf{C2%
}^{\prime }.$ The following assertions are equivalent:

\begin{enumerate}
\item Every subtractive ideal of $S$ is a direct summand;

\item Every $S$-semimodule is $S$-$e$-injective;

\item Every $S$-semimodule is $S$-$i$-injective;

\item Every subtractive ideal of $S$ is $S$-$i$-injective;

\item Every short exact sequence of $S$-semimodules $0\rightarrow
L\rightarrow S\rightarrow N\rightarrow 0$ is left splitting;

\item $S$ is $k$-Noetherian;

\item $S$ satisfies ACC on direct summands;

\item $S$ satisfies DCC on direct summands;

\item $S=S_{1}\oplus S_{2}\oplus S_{3}\oplus ...\oplus S_{n}$, where every $%
S_{i}$ is an irreducible summand;

\item $S$ is congruence-semisimple.
\end{enumerate}
\end{thm}

Combining Theorems \ref{cong-C2} and \ref{CSS-COMM}, we obtain the following
characterization of commutative congruence-semisimple semirings:

\begin{thm}
\label{com-css}The following assertions are equivalent for a \emph{%
commutative} semiring $S$ which satisfies $\mathbf{C2}^{\prime }$:

\begin{enumerate}
\item Every subtractive ideal of $S$ is a direct summand;

\item Every $S$-semimodule is $S$-$e$-projective ($S$-$k$-projective);

\item Every $S$-semimodule is $S$-$e$-injective ($S$-$i$-injective);

\item For every subtractive ideal $I$ of $S$ we have: $S/I$ is $S$-$k$%
-projective ($I$ is $S$-$i$-injective);

\item Every short exact sequence $0\longrightarrow I\longrightarrow
S\longrightarrow N\longrightarrow 0$ in $_{S}\mathbf{SM}$ right splits (left
splits);

\item $S$ is $k$-Noetherian;

\item $_{S}S$ satisfies the ACC on the direct summands;

\item $_{S}S$ satisfies the DCC on the direct summands;

\item $S=S_{1}\oplus S_{2}\oplus S_{3}\oplus ...\oplus S_{n}$, where every $%
S_{i}$ is an irreducible summand;

\item $S$ is congruence-semisimple.
\end{enumerate}
\end{thm}

The following example shows that the assumption that $S$ is a \emph{%
commutative }semiring cannot be dropped from Theorems \ref{id-SS-C1}, \ref%
{ISS-COMM}, whence not from our main result\ Theorems \ref{com-idss}, \ref%
{com-css}:

\begin{ex}
\label{nonseproj}Consider the semiring $S:=M_{2}(\mathbb{R}^{+})$.

\begin{enumerate}
\item $S$ is a left ideal-semisimple semiring.

\item $N_{1}$ is a subtractive left ideal of $S$ which is not a direct
summand.

\item $S/N_{1}$ is not an $S$-$k$-projective $S$-semimodule (whence not $S$-$%
e$-projective).

\item $N_{1}$ is not $S$-$e$-injective.
\end{enumerate}
\end{ex}

\begin{Beweis}
\begin{enumerate}
\item The semiring $M_{2}(\mathbb{R}^{+})$ is left ideal-simple since $%
\mathbb{R}^{+}$ is a semifield (\cite[Theorem 7.8]{HW1998}).

\item Let $K$ be a left ideal of $S$ such that $S=N_{1}+K$. Then $1_{S}=i+k$
for some $i\in N_{1}$ and $k\in K$, that is
\begin{equation*}
\left[ {%
\begin{array}{cc}
1 & 0 \\
0 & 1%
\end{array}%
}\right] =\left[ {%
\begin{array}{cc}
a & a \\
b & b%
\end{array}%
}\right] +\left[ {%
\begin{array}{cc}
p & q \\
r & s%
\end{array}%
}\right] .
\end{equation*}%
Then $p+a=1=s+b$ and $q+a=0=r+b$, whence $a=q=r=b=0$ as $\mathbb{R}^{+}$ is
zerosumfree. Therefore, $i=0$ and $k=1_{S}$, which implies $K=S$ and $0\neq
N_{1}=N_{1}\cap K.$ Thus, the sum $N_{1}+K$ is not direct. Consequently, $%
N_{1}$ is a subtractive left ideal of $S$ which is not a direct summand.

\item Let $\pi :S\rightarrow S/N_{1}$ be the canonical map and $id_{S/N_{1}}$
be the identity map of $S/N_{1}$. Notice that $\pi $ is a normal
epimorphism. Consider
\begin{equation*}
e_{1}=\left[ {%
\begin{array}{cc}
1 & 0 \\
0 & 0%
\end{array}%
}\right] \text{ and }e_{2}=\left[ {%
\begin{array}{cc}
0 & 0 \\
0 & 1%
\end{array}%
}\right] .
\end{equation*}%
Suppose that there exists an $S$-linear map $g:S/N_{1}\rightarrow S$ such
that $\pi \circ g=id_{S/N_{1}}$. Then $g(\overline{e_{1}})\in \pi ^{-1}(%
\overline{e_{1}})$ and $g(\overline{e_{2}})\in \pi ^{-1}(\overline{e_{2}})$.
Write $g(\overline{e_{1}}):=\left[ {%
\begin{array}{cc}
p & q \\
r & s%
\end{array}%
}\right] $ for some $p,q,r,s\in \mathbb{R}^{+}$. Then $\left[ {%
\begin{array}{cc}
p+k & q+k \\
r+l & s+l%
\end{array}%
}\right] =\left[ {%
\begin{array}{cc}
m+1 & m \\
n & n%
\end{array}%
}\right] $ for some $k,l,m,n\in \mathbb{R}^{+}$, whence $r=s$ and $p=q+1$ as
$\mathbb{R}^{+}$ is cancellative. By relabeling, we have $g(\overline{e_{1}}%
)=\left[ {%
\begin{array}{cc}
a+1 & a \\
b & b%
\end{array}%
}\right] $ for some $a,b\in \mathbb{R}^{+}$. Similarly, $g(\overline{e_{2}})=%
\left[ {%
\begin{array}{cc}
c & c \\
d & d+1%
\end{array}%
}\right] $ for some $c,d\in \mathbb{R}^{+}$.

Let $x:=\left[ {%
\begin{array}{cc}
p & q \\
r & s%
\end{array}%
}\right] \in S.$ Then $x=\left[ {%
\begin{array}{cc}
p & 0 \\
r & 0%
\end{array}%
}\right] e_{1}+\left[ {%
\begin{array}{cc}
0 & q \\
0 & s%
\end{array}%
}\right] e_{2},$ whence%
\begin{equation*}
g(\overline{x})=\left[ {%
\begin{array}{cc}
p & 0 \\
r & 0%
\end{array}%
}\right] g(\overline{e_{1}})+\left[ {%
\begin{array}{cc}
0 & q \\
0 & s%
\end{array}%
}\right] g(\overline{e_{2}})=\left[ {%
\begin{array}{cc}
pa+dq+p & pa+dq+q \\
ra+sd+r & ra+sd+s%
\end{array}%
}\right] .
\end{equation*}%
But $x=\left[ {%
\begin{array}{cc}
p & 1 \\
r & 0%
\end{array}%
}\right] e_{1}+\left[ {%
\begin{array}{cc}
0 & q \\
1 & s%
\end{array}%
}\right] e_{2}$, whence%
\begin{equation*}
\begin{array}{ccc}
\left[ {%
\begin{array}{cc}
pa+dq+p & pa+dq+q \\
ra+sd+r & ra+sd+s%
\end{array}%
}\right] & = & g(\overline{x})=\left[ {%
\begin{array}{cc}
p & 1 \\
r & 0%
\end{array}%
}\right] g(\overline{e_{1}})+\left[ {%
\begin{array}{cc}
0 & q \\
1 & s%
\end{array}%
}\right] g(\overline{e_{2}}) \\
& = & \left[ {%
\begin{array}{cc}
(pa+dq+p)+b & (pa+dq+q)+b \\
(ra+sd+r)+c & (ra+sd+s)+c%
\end{array}%
}\right] ,%
\end{array}%
\end{equation*}%
whence $b=0=c$ as $\mathbb{R}^{+}$ is cancellative. Thus $g(\overline{e_{1}}%
)=\left[ {%
\begin{array}{cc}
a+1 & a \\
0 & 0%
\end{array}%
}\right] $ for some $a,b\in \mathbb{R}^{+}$ and $g(\overline{e_{2}})=\left[ {%
\begin{array}{cc}
0 & 0 \\
d & d+1%
\end{array}%
}\right] $.

Let $y=\left[ {%
\begin{array}{cc}
2 & 1 \\
0 & 0%
\end{array}%
}\right] $. Notice that $\overline{e_{1}}=\overline{y}$, whence
\begin{equation*}
\left[ {%
\begin{array}{cc}
a+1 & a \\
0 & 0%
\end{array}%
}\right] =g(\overline{e_{1}})=g(\overline{y})=\left[ {%
\begin{array}{cc}
2a+d+2 & 2a+d+1 \\
0 & 0%
\end{array}%
}\right] ,
\end{equation*}%
and so $a=2a+d+1$. Since $\mathbb{R}^{+}$ is cancellative, $a+d+1=0$, that
is $1$ has an additive inverse, a contradiction. Hence, there is no such $S$%
-linear map $g$ with $\pi \circ g=id_{S/I};$ \emph{i.e.}, $S/I$ is not $S$-$%
k $-projective. Since $S/I$ is not $S$-$k$-projective, $S/I$ is not $S$-$e$%
-projective.

\item This was shown in \cite[Example 2.19.]{AN-2}.$\blacksquare $
\end{enumerate}
\end{Beweis}

%\addcontentsline{toc}{section}{\protect\numberline{}{Index}} \printindex

\end{document}